\documentclass[10pt]{amsart}

\usepackage{amsmath,amssymb,epsfig,leftidx}
\usepackage{float, color}
\usepackage{amsmath}
\usepackage{amssymb}
\usepackage{graphicx}
\usepackage{graphics}

\theoremstyle{definition}

\theoremstyle{remark}

\begin{document}

\title[Cone beam vector tomography]{An improved exact inversion formula for cone beam vector tomography}
\author{Alexander Katsevich, Dimitri Rothermel and Thomas Schuster}
\address{Alexander Katsevich: Department of Mathematics, University of Central Florida, Orlando, FL 32816-1364;}
\address{Dimitri Rothermel, Thomas Schuster: Department of Mathematics, Saarland University, 66123 Saarbr\"ucken, Germany}

\email{Alexander.Katsevich@ucf.edu}
\email{thomas.schuster@num.uni-sb.de}
\email{rothermel@math.uni-sb.de}

\numberwithin{equation}{section}

\newcommand{\fbfeins}{\mathbf{f}^{(1)}}
\newcommand{\fbfzwei}{\mathbf{f}^{(2)}}
\newcommand{\ang}{\theta}
\newcommand{\vct}{\Theta}
\newcommand{\vctp}{\vct^{\perp}}
\newcommand{\sw}{\text{supp}\,w}
\newcommand{\xo}{x^{(0)}}
\newcommand{\xio}{\xi^{(0)}}
\newcommand{\zo}{z^{(0)}}
\newcommand{\xs}{x^{(s)}}
\newcommand{\rtc}{\br^3\setminus C}
\newcommand{\cD}{{\mathcal D}}
\newcommand{\cS}{{\mathcal S}}
\newcommand{\cBf}{{\mathcal B}_{\eta}f}
\newcommand{\CA}{{\mathcal A}}
\newcommand{\CB}{{\mathcal B}}
\newcommand{\xyf}{\frac{x-y(s)}{|x-y(s)|}}
\newcommand{\cx}{\check x}
\newcommand{\cy}{\check y}
\newcommand{\cz}{\check z}
\newcommand{\cxi}{\check \xi}
\newcommand{\ap}{\al,p}
\newcommand{\fbf}{\mathbf{f}}
\newcommand{\bp}{\boldsymbol\phi}
\newcommand{\geg}{C_n^{(3/2)}}
\newcommand{\Dbf}{\mathbf{D}}
\newcommand{\drm}{\mathrm{d}}
\newcommand{\RR}{\mathbb{R}}
\newcommand{\Rbf}{\mathbf{R}}
\newcommand{\Rtan}{\mathbf{R}^{(tan)}}
\newcommand{\Rnor}{\mathbf{R}^{(nor)}}
\newcommand{\Pbf}{\mathbf{P}}
\newcommand{\al}{\alpha} 
\newcommand{\bt}{\beta} 
\newcommand{\de}{\delta} 
\newcommand{\ga}{\gamma} 
\newcommand{\G}{\Gamma} 
\newcommand{\e}{\epsilon} 
\newcommand{\la}{\lambda} 
\newcommand{\si}{\sigma} 
\newcommand{\ph}{\varphi} 
\newcommand{\pa}{\partial}       
\newcommand{\vn}{\varnothing}       
\newcommand{\coi}{C_0^{\infty}} 
\newcommand{\ci}{C^{\infty}} 
\newcommand{\ts}{\text{supp}\,} 
\newcommand{\tss}{\text{singsupp}\,} 
\newcommand{\br}{{\mathbb R}} 
\newcommand{\rn}{\mathbb R^n} 
\newcommand{\adx}{\vct\cdot x} 
\newcommand{\Sn}{S^{n-1}} 
\newcommand{\sg}{\text{sgn}} 
\newcommand{\ibr}{\int_{\br}} 
\newcommand{\brt}{\br^2}       
\newcommand{\irt}{\int_{\brt}}       
\newcommand{\ioi}{\int_0^{\infty}}       
\newcommand{\iii}{\int_{-\infty}^{\infty}}       
\newcommand{\isn}{\int_{\Sn}}       
\newcommand{\fle}{f_{\Lambda\e}}
\newcommand{\CS}{{\mathcal S}}
\newcommand{\CD}{{\mathcal D}}
\newcommand{\bphi}{{\boldsymbol{\Phi}}}
\newcommand{\bpsi}{{\boldsymbol{\psi}}}

\begin{abstract}
In this article we present an improved exact inversion formula for the 3D cone beam transform of vector fields.
It is well known that only the solenoidal part of a vector field can be determined by the longitudinal ray transform
of a vector field in cone beam geometry.
The exact inversion formula, as it was developed in A. Katsevich and T. Schuster, {\it An exact inversion formula for cone
beam vector tomography}, Inverse Problems 29 (2013), consists of two parts. the first part is of filtered
backprojection type, whereas the second part is a costly 4D integration and very inefficient. In this article
we tackle this second term and achieve an improvement which is easily to implement and saves one order of
integration. The theory says that the first part contains all information about the curl of the field, whereas
the second part presumably has information about the boundary values. This suggestion is supported by the fact
that the second part vanishes if the exact field is divergence free and tangential at the boundary. A number
of numerical tests, that are also subject of this article, confirm the theoretical results and the exactness of
the formula.
\end{abstract}

\keywords{Cone beam, vector tomography, theoretically exact reconstruction, general trajectory}

\subjclass{Primary 44A12, 65R10, 92C55}
\maketitle

\section{Introduction}\label{intro}

We consider the problem of reconstructing a smooth vector field $\fbf$, supported in the open unit ball $B^3=\{ x\in\RR^3 : |x| < 1\} \subset \RR^3$,
from its cone beam data
\begin{equation}\label{CBT}
[\Dbf \fbf]\big(y(s),\Theta\big) = \int\limits_0^\infty \fbf \big(y(s)+t \Theta\big)\cdot\Theta\,d t,\ \Theta\in S^2.
\end{equation}
Here $y(s)$, $s\in \Lambda\subset \RR$, denotes a parametrization of the source
trajectory $\Gamma\subset (\RR^3\backslash {\overline {B^3}})$ and $S^2:=\partial B^3$ is the unit sphere in $\br^3$. It is assumed that $\Theta\in \mathcal{C}$, where $\mathcal{C}$ is a cone
and that $B^3 \subset y(s)+\mathcal{C}$, i.e. the unit ball is completely
contained inside the union of the rays emanating from any source position $y(s)$. The cone beam transform \eqref{CBT} is the mathematical model
of vector tomography, where a flow field $\fbf$ is reconstructed from ultrasound Doppler or time-of-flight measurements with sources located on the trajectory $\Gamma$,
see e.g. \cite{SPARR;STRAHLEN;LINDSTROEM;PERSSON:1995,SCHUSTER;THEIS;LOUIS:2009}. It is well known that $\Dbf$
has a non-trivial null space and that only the solenoidal part $\fbf^s$ of $\fbf$ can be reconstructed
from $\Dbf \fbf$ \cite{shara_94}. 

In \cite{katsch13} the authors obtained the first explicit and theoretically exact inversion formula for the cone beam transform of vector fields, which is not based on series expansions. The formula gives an analytical expression for computing the solenoidal part $\fbf^s$ of $\fbf$ from $\Dbf \fbf$. The inversion formula consists of two parts: $\mathbf{f}^s=\fbf_1+\fbf_2$. The first part that recovers $\fbf_1$ is of convolution-backprojection type,
\begin{equation}
\begin{split}
{\bf f}_1(x)
=\frac1{8\pi^2}\int_I & \frac1{|x-y(s)|} 
\int_0^{2\pi} [\bphi_{\theta\theta}(s,\al(\theta))+\bphi(s,\al(\theta))] \\
& \times\int_0^{2\pi}  \frac{g(y(s),\cos\ga\,\al^\perp(\theta) + \sin\ga\,\bt) }{\cos^2\ga} d\ga d\theta
ds,
\end{split}
\label{recon-3}
\end{equation}
with $\beta=\beta(s,x)=(x-y(s))/|x-y(s)|$ and $\bphi$ can be computed from the data ${\mathbf D}\fbf$. The second part that computes $\fbf_2$ is much more complex and less efficient. It consists of a costly 4D integral over $S^2\times S^2$ and resembles the early approaches to inverting the cone beam transform based on the Tuy and Grangeat formulas \cite{GRANGEAT:1991, ks1, ZENG;CLACK;GULLBERG:94}. The main result of this paper is the development of an efficient formula for computing $\fbf_2$.  

The paper is organized as follows. In Section~\ref{derivation} we obtain a new formula for computing $\fbf_2$ in the case when the support of $\mathbf{f}$ is the unit ball. Then, in Section~\ref{general_dom}, we outline an algorithm for computing $\fbf_2$ for general domains. The results of numerical testing of the formula for $\fbf_2$ are presented in Section~\ref{numerics}. Testing of the algorithm for computing $\fbf_2$ for general domains will be the subject of future research.

\section{Derivation}\label{derivation}

The Radon transform of ${\bf f}^s$ can be written in the form \cite{ks-11}
\begin{equation}
\begin{split}
\Rbf{\bf f}^s(p,\eta)=&8\pi^2\sum_{n\ge0} \frac{ (1-p^2)\geg(p)}{(n+1)(2n+3)} 
\sum_{|l|\le n+1} b_{n+1,l}^{(n)} ((n+1){\bf y}^{(1)}_{n+1,l}(\eta)+{\bf y}^{(2)}_{n+1,l}(\eta))\\
&+(\text{${\bf y}^{(2)}$ terms with different indices and ${\bf y}^{(3)}$ terms}).
\end{split}
\label{rfexp}
\end{equation}
Here $\geg$ are the Gegenbauer polynomials, and ${\bf y}^{(j)}_{n,l}$, $j=1,2,3$, are the vector spherical harmonics (see \cite{dks-07, ks-11}). Differentiating \eqref{rfexp} with respect to $p$ and using the identity:
\begin{equation}
[(1-p^2)\geg(p)]''=-(n+1)(n+2)\geg(p),
\label{geg-der}
\end{equation}
the second derivative of the Radon transform of ${\bf f}^s$ is given by
\begin{equation}
\begin{split}
\pa_p^2 \Rbf{\bf f}^s(p,\eta)=&-8\pi^2\sum_{n\ge0} \frac{ (n+2)\geg(p)}{2n+3} 
\sum_{|l|\le n+1} b_{n+1,l}^{(n)} ((n+1){\bf y}^{(1)}_{n+1,l}(\eta)+{\bf y}^{(2)}_{n+1,l}(\eta))\\
&+(\text{${\bf y}^{(2)}$ terms with different indices and ${\bf y}^{(3)}$ terms}).
\end{split}
\label{rfexp_v2}
\end{equation}

Pick a ``reasonable function" $\phi(p)$ defined on $[-1,1]$, multiply (\ref{rfexp}) by $\phi(p){\overline {\bf y}^{(2)}_{n+1,l}}(\eta)$, and integrate over $[-1,1]\times S^2$. Here and below the overbar denotes complex conjugation. Since vector spherical harmonics are orthogonal, we get
\begin{equation}
\begin{split}
\int_{S^2} \int_{-1}^1 & \pa_p^2\Rbf{\bf f}^s(p,\eta)\cdot [\phi(p){\overline {\bf y}^{(2)}_{n+1,l}}(\eta)] dp d\eta 
= -\frac{8\pi^2 \check\phi_n (n+2) b_{n+1,l}^{(n)}}{2n+3}  \lVert{\bf y}^{(2)}_{n+1,l}\rVert^2,\\
\check\phi_n:&=   \int_{-1}^1 \geg(p)\phi(p)dp.
\end{split}
\label{coefs}
\end{equation}

Since ${\bf y}^{(2)}_{n+1,l}(\eta)$ is orthogonal to the normal component $\Rnor \fbf$ of $\Rbf{\bf f}$ (see \cite{ks-11}), we can replace $\Rbf{\bf f}^s$ with $\Rtan{\bf f}^s$ in (\ref{coefs}). The latter is equal to $\Rtan{\bf f}$, as follows from the orthogonal expansions used in \cite{ks-11}. Here $\Rnor \fbf$ and $\Rtan{\bf f}$ are the normal and tangential components of the Radon transform of $\fbf$, respectively (cf. \cite{ks-11, katsch13}):
\begin{equation}
\begin{split}
[\Rnor \fbf](p,\eta) := (\eta\cdot [\Rbf\fbf ](p,\eta))\eta,\
[\Rtan \fbf](p,\eta) := [\Rbf\fbf](p,\eta) - [\Rnor\fbf](p,\eta).
\end{split}
\end{equation}
Multiply the top equation in (\ref{coefs}) by 
\begin{equation}
(n+1)\frac{\geg(q)}{\check\phi_n}\frac{{\bf y}^{(1)}_{n+1,l}(\al)}{\lVert{\bf y}^{(2)}_{n+1,l}\rVert^2}
\label{mult-factor}
\end{equation}
and sum over $n,l$ to get the derivative $\partial_p^2\Rnor{\bf f}^s$:
\begin{equation}
\begin{split}
[\partial_q^2\Rnor{\bf f}^s](q,\al)
=\sum_{n\ge0} \sum_{|l|\le n+1} &
\int_{S^2} \int_{-1}^1 \pa_p^2\Rtan
{\bf f}(p,\eta)\cdot [\phi(p){\overline {\bf y}^{(2)}_{n+1,l}}(\eta)] dp d\eta \\
&\times(n+1)\frac{\geg(q)}{\check\phi_n}\frac{{\bf y}^{(1)}_{n+1,l}(\al)}{\lVert{\bf y}^{(2)}_{n+1,l}\rVert^2}.
\end{split}
\label{rnorm}
\end{equation}
Similarly to \cite{katsch13},  using that ${\bf y}^{(1)}_{n,l}(\al)=\al Y_{n,l}(\al)$, ${\bf y}^{(2)}_{n,l}(\eta) =\nabla_\eta Y_{n,l}(\eta)$, and $\lVert{\bf y}^{(2)}_{n+1,l}\rVert^2=(n+1)(n+2)$ (see \cite{dks-07, ks-11}), we write the sum with respect to $l$ in the form of a rank-one matrix
\begin{equation}
\al\otimes\nabla_\eta \sum_{|l|\le n+1} Y_{n+1,l}(\al) {\overline Y}_{n+1,l}(\eta)
=\frac{2n+3}{4\pi}\al\otimes\nabla_\eta P_{n+1}(\al\cdot\eta).
\label{lsum}
\end{equation}
Here $P_n$ are the Legendre polynomials, $Y_{n,l}$ are the scalar spherical harmonics, and we used the addition theorem for spherical harmonics. The operator in (\ref{lsum}) acts on vectors by computing the dot product of an input vector with $\nabla_\eta P_{n+1}(\al\cdot\eta)$ and then multiplying the result by the vector $(2n+3)\al/4\pi$. 

Using (\ref{lsum}) in (\ref{rnorm}), substituting the result into the Radon transform inversion formula, and (so far formally) changing the order of integration and summation  we get
\begin{equation}
\begin{split}
&\fbf^{(2)} (x) = -\frac{1}{4(2\pi)^3} 
\int_{S^2}  
K(x,\eta) \int_{-1}^1 \phi(p)\pa_p^2\Rtan
{\bf f}(p,\eta)
dp d\eta,\\
&K_1(x,\alpha;\eta):=\sum_{n\ge0} 
\frac{2n+3}{\check\phi_n(n+2)} P_{n+1}(\al\cdot\eta)\geg(\al\cdot x),\\
&K(x,\eta):=\int_{S^2} \al\otimes \nabla_\eta K_1(x,\alpha;\eta)d\al.
\end{split}
\label{rnorm-1}
\end{equation} 
Define 
\begin{equation}
\phi(p):=(1-p^2)\sum_{n\ge 0}\frac{\check\phi_n}{\Vert \geg\Vert^2} \geg(p).
\label{pfd2}
\end{equation} 
From the orthogonality of the Gegenbauer polynomials (eq. 22.2.3 in \cite{as}), $\phi(p)$ satisfies (see the second equation in \eqref{coefs})
\begin{equation}
\check\phi_n=\int_{-1}^1 \geg(p)\phi(p)dp.
\label{coefs2}
\end{equation}
In view of \eqref{coefs2}, any function given by \eqref{pfd2} can be used in \eqref{rnorm-1}. 
The goal is to choose the coefficients $\check\phi_n$ so we could use the following identity (see 5.10.2.2 in \cite{pbm2}), which we write here in a symmetric form:
\begin{equation}
\begin{split}
\sum_{n\ge1}\frac{2n+1}{n(n+1)}P_n(s) P_n(t)&=2\ln2-1-\ln(1+|t-s|-st)\\
& =2\ln2-1-\ln((1+\max(s,t))(1-\min(s,t))),
\end{split}
\label{PPid}
\end{equation}
$-1 <s, t <1$. In view of the relation (eq. 22.5.37 in \cite{as})
\begin{equation}
P_{n+1}'(t)=\geg(t),
\label{Pgeg}
\end{equation}
we have to evaluate the expression
\begin{equation}
S_2(s,t):=\sum_{n\ge1} 
\frac{2n+1}{\check\phi_{n-1}(n+1)} P_n(s)P_n(t).
\label{S2kern}
\end{equation} 
Comparing \eqref{PPid} and \eqref{S2kern} we select $\check\phi_n=n+1$. 

Substituting into \eqref{pfd2} and using that $\Vert \geg\Vert^2=(n+1)(n+2)/(n+(3/2))$ gives
\begin{equation}
\phi(p)=\frac{1-p^2}2 \sum_{n\ge 0}\frac{2n+3}{n+2} \geg(p)
=\frac{1-p^2}2 \sum_{n\ge 0}\frac{2n+3}{n+2} P_{n+1}'(p).
\label{pfdser}
\end{equation} 
Consider the integral with respect to $p$ in \eqref{rnorm-1}. Define ${\bf u}(p):=\pa_p^2\Rtan {\bf f}(p,\eta)$. By assumption, ${\bf u}$ is smooth in $p$. Ignoring the dependence on $\eta$ we can write this integral in the form
\begin{equation}
\int_{-1}^{1} \phi(p){\bf u}(p)dp=\lim_{N\to\infty} \frac12 \int_{-1}^{1} \left[\sum_{n= 0}^N\frac{2n+3}{n+2} P_{n+1}(p)\right]((1-p^2){\bf u}(p))'dp.
\label{p-int}
\end{equation} 
Surprisingly, the expression in brackets is exactly the same as the one occuring in (3.26), (3.29) of \cite{katsch13}. The derivation (3.28)--(3.41) of \cite{katsch13} justifies taking the limit inside the integral in \eqref{p-int}. Using (3.42) of \cite{katsch13} and ignoring the constant terms because of the derivative in \eqref{p-int} gives:
\begin{equation}
\begin{split}
\int_{-1}^{1} & \phi(p){\bf u}(p)dp\\
&=\frac12 \int_{-1}^{1} \left[\sqrt{2/(1-p)}-\ln(1+\sqrt{(1-p)/2})+\frac12\ln(1-p)\right]((1-p^2){\bf u}(p))'dp.
\end{split}
\label{p-int-expl}
\end{equation} 
Integrating by parts again we immediately get:
\begin{equation}
\begin{split}
\phi(p)
=(1-p^2)\left[(2-2p)^{-3/2}-\frac12\left((2-2p)^{-1/2}(1-p+\sqrt{2-2p}\,)^{-1}\right)\right].
\end{split}
\label{twosers}
\end{equation} 
Clearly, $\phi\in L^1([-1,1])$.


Now we can find the kernel $K$. Denote 
\begin{equation}
{\bf v}(\eta):=\int_{-1}^1 \phi(p)\pa_p^2\Rtan{\bf f}(p,\eta)dp.
\label{psi-def}
\end{equation} 
By assumption, ${\bf v}\in\ci(S^2)$. Suppose first that ${\bf v}$ is a linear combination of the first $N_0$ vector spherical harmonics. Using \eqref{Pgeg}, we can rewrite the integral in the first line of \eqref{rnorm-1} as follows:
\begin{equation}
\begin{split}
&{\bf V}(x):= \nabla_x v(x),\ v(x):=\int_{S^2}  \int_{S^2}
\nabla_\eta \tilde K_N(x,\al;\eta)\cdot{\bf v}(\eta) d\eta d\al,\\
&\tilde K_N(x,\alpha;\eta):=\sum_{n=1}^N \frac{2n+1}{n(n+1)} P_n(\al\cdot\eta)P_n(\al\cdot x),
\end{split}
\label{step1}
\end{equation} 
for any $N\ge N_0$. Integrating by parts on the unit sphere gives (see e.g. (3.27) in \cite{katsch13}):
\begin{equation}
v(x):=\int_{S^2}  \int_{S^2}
\tilde K_N(x,\al;\eta)(L_\eta{\bf v})(\eta) d\eta d\al,
\label{step2}
\end{equation} 
with the differential operator $(L_\eta {\bf v})(\eta) := (2\eta-\nabla_\eta)\cdot {\bf v}(\eta)$. Given that $|P_n(t)|\le 1$ and $|P_n(t)|=O(n^{-1/2})$ uniformly on compact subsets of $(-1,1)$ (which follows from the inequality 8.917.4 in \cite{gr}), we can take the limit as $N\to\infty$ inside the integral because the series is absolutely convergent as long as $|x|<1$. Hence \eqref{PPid} implies
\begin{equation}
v(x)=\int_{S^2}  \int_{S^2} S_2(\al\cdot x,\al\cdot \eta)(L_\eta{\bf v})(\eta) d\eta d\al,\ 
S_2(s,t):=-\ln(1+|t-s|-ts),
\label{step3}
\end{equation} 
where we again ignored the constant terms because of the derivatives in \eqref{step3}.
Integration by parts in the sense of distributions converts $L_\eta$ back into $\nabla_\eta$. Since $\nabla_\eta=\nabla-(\eta\cdot\nabla)\eta$, it is easy to check that for a differentiable function $F$ defined on $\br$ we have $\nabla_\eta F(\al\cdot\eta)=F'(\al\cdot\eta)(\al-(\al\cdot\eta)\eta)$. Thus, 
\begin{equation}
\label{step4}
v(x) = \int_{S^2}  \int_{S^2} \left.\frac{\pa}{\pa t}S_2(\al\cdot x,t)\right|_{t=\al\cdot\eta}(\al-(\al\cdot\eta)\eta)\cdot{\bf v}(\eta) d\eta d\al,
\end{equation}
where
\[
\frac{\pa}{\pa t}S_2(s,t)=-\frac{\pa}{\pa t}\ln(1+|t-s|-ts)
=-\frac{\text{sgn}(t-s)-s}{1+|t-s|-ts}
=-\frac{1}{t+\text{sgn}(t-s)}.
\]

Equation \eqref{step4} implies that $v(x)$ is the result of applying a distribution, which depends smoothly on the parameter $x$, to the test function $(\al-(\al\cdot\eta)\eta)\cdot{\bf v}(\eta) \in\ci(S^2\times S^2)$. Substitute $s=\al\cdot x$ into the formula for $\pa S_2/\pa t$ in \eqref{step4}. An easy calculation shows that in the sense of distributions:
\begin{equation}
\nabla_x\frac{-1}{t+\text{sgn}(t-\al\cdot x)}
=\al\left.\frac{\pa}{\pa s}\frac{-1}{t+\text{sgn}(t-s)}\right|_{s=\al\cdot x}
=\al \frac{2\de(\al\cdot x-t)}{1-t^2}.
\label{step_aux}
\end{equation} 
Combining  \eqref{step2}--\eqref{step_aux} and  \eqref{rnorm-1} we obtain that the operator $K$ is given by 
\begin{equation}
K(x,\eta)=2\int_{S^2} \al\otimes[\al-(\al\cdot\eta)\eta] \frac{\de(\al\cdot(\eta-x))}{1-(\al\cdot x)^2}d\al.
\label{K-kern}
\end{equation} 
Thus, the formal calculation in \eqref{rnorm-1} is justified. For convenience, we replaced $\al\cdot \eta$ with $\al\cdot x$ in \eqref{K-kern}. The two are equal on the support of the delta function.

Next we compute the kernel $K$ explicitly. As is easily seen, $K(x,\eta)\eta\equiv 0$. Introduce the coordinate system in which $\eta=e_3$, and $x$ lies in the $e_1,e_3$-plane. This implies that the matrix $K(x,\eta)$ has the following zero components: $K_{13}=K_{23}=K_{33}=0$. 

The integral in \eqref{K-kern} is over the great circle orthogonal to $\eta-x$. Consider two points on that circle with the same first coordinates. Clearly, the third coordinates of these two points will also be equal to each other, but they will have opposite (i.e., equal in magnitude and of opposite signs) second coordinate. This implies that $K_{12}=K_{21}=K_{32}=0$. 

To compute the remaining nonzero components  $K_{11},K_{22},K_{31}$, let $\mu$ denote the angle between the vectors $\eta$ and $\eta-x$.  Let $\theta$ denote the polar angle in the plane orthogonal to $\eta-x$. Denote also $L:=|\eta-x|$. The points on the great circle are parameterized as follows:
\begin{equation}
(\cos\mu\cos\theta,\sin\theta,\sin\mu\cos\theta),\ 0\le\theta\le 2\pi.
\label{grc}
\end{equation} 
Note that the term $(\al\cdot\eta)\eta$ in the numerator in \eqref{K-kern} does not contribute to the components we need to calculate. Therefore, using the homogeneity of the delta-function and the following formulas
\begin{equation}
\begin{split}
\frac{\cos^2\theta}{1-r^2\cos^2\theta}
&=\frac{1}{r^2}\left(-1+\frac1{1-r^2\cos^2\theta}\right),\\
\frac{\sin^2\theta}{1-r^2\cos^2\theta}
&=\frac{1}{r^2}\left(1-\frac{1-r^2}{1-r^2\cos^2\theta}\right),\\
\int_{-\pi}^{\pi}\frac{d\theta}{1-r^2\cos^2\theta}&=\frac{2\pi}{\sqrt{1-r^2}},
\end{split}
\label{simple-ids}
\end{equation} 
we find with $r=\sin\mu$
\begin{equation}
\begin{split}
LK_{11}(x,\eta)&=\cos^2\mu\int_0^{2\pi} \frac{\cos^2\theta}{1-\sin^2\mu\cos^2\theta}d\theta
=\frac{\cos^2\mu}{\sin^2\mu}\left(-2\pi+\frac{2\pi}{\sqrt{1-\sin^2\mu}}\right)\\
&=2\pi \frac{\cos\mu(1-\cos\mu)}{\sin^2\mu}=2\pi \frac{\cos\mu}{1+\cos\mu}.
\end{split}
\label{K11}
\end{equation} 
In a similar fashion,
\begin{equation}
\begin{split}
LK_{22}(x,\eta)&=\int_0^{2\pi} \frac{\sin^2\theta}{1-\sin^2\mu\cos^2\theta}d\theta
=\frac{1}{\sin^2\mu}\left(2\pi-\frac{2\pi(1-\sin^2\mu)}{\sqrt{1-\sin^2\mu}}\right)\\
&=2\pi \frac{1-\cos\mu}{\sin^2\mu}=2\pi \frac{1}{1+\cos(\mu)},
\end{split}
\label{K22}
\end{equation} 
and
\begin{equation}
\begin{split}
LK_{31}(x,\eta)&=\sin\mu\cos\mu\int_0^{2\pi} \frac{\cos^2\theta}{1-\sin^2\mu\cos^2\theta}d\theta
=2\pi \frac{\sin\mu}{1+\cos\mu}.
\end{split}
\label{K31}
\end{equation} 
Application of the matrix $K$ to a vector $h$ is given by
\begin{equation}
\begin{split}
K(x,\eta)h=\frac{2\pi}{L(1+\cos\mu)}((h\cdot e_1)\cos\mu e_1+(h\cdot e_2)e_2+(h\cdot e_1)\sin\mu e_3)
\end{split}
\label{Kh}
\end{equation} 
As is easily checked, $u:=\cos\mu e_1+\sin\mu e_3$ is the unit vector perpendicular to $\eta-x$ and lying in the $\eta,x$-plane. Therefore,
\begin{equation}
\begin{split}
K(x,\eta)h=\frac{2\pi}{L(1+\cos\mu)}((h\cdot e_1)u+(h\cdot e_2)e_2).
\end{split}
\label{Kh2}
\end{equation} 
In coordinate-free form we have
\begin{equation}
\begin{split}
e_1=\frac{x-(\eta\cdot x)\eta}{|\eta \times x|},\
e_2=\frac{\eta\times x}{|\eta\times x|},\
e_3=\eta,\ u=e_v\times \frac{\eta\times x}{|\eta\times x|},\ e_v:=\frac{\eta-x}{|\eta-x|}.
\end{split}
\label{vectors}
\end{equation} 
Also,
\begin{equation}
L(1+\cos\mu)=|\eta-x|+\eta\cdot(\eta-x),\ L\sin\mu=|\eta \times x|.
\label{denom}
\end{equation} 
Since $u,e_v$, and $e_2$ form an orthonormal triple, 
\begin{equation}
(h\cdot e_1)u+(h\cdot e_2)e_2=h-(h\cdot e_v)e_v-(h\cdot (u-e_1))u.
\label{hexp}
\end{equation}

Combining \eqref{rnorm-1}, \eqref{K-kern}, \eqref{Kh2}, \eqref{vectors}, \eqref{denom}, the formula for $\fbf^{(2)}$ becomes
\begin{equation}
\begin{split}
\fbf^{(2)} (x) &= -\frac{1}{8\pi^2} 
\int_{S^2}  \frac{(\Psi(\eta)\cdot e_1)u+(\Psi(\eta)\cdot e_2)e_2}{|\eta-x|+\eta\cdot(\eta-x)} d\eta,\\
\Psi(\eta):&= \int_{-1}^1 \phi(p)\pa_p^2\Rtan{\bf f}(p,\eta)dp.
\end{split}
\label{f2-final}
\end{equation} 

A disadvantage of the formula \eqref{f2-final} is that it appears to have non-smooth dependence on $x$ in a neighborhood of $x=0$. Indeed, if $x=0$, then a number of vectors in \eqref{vectors} are undefined. Thus, we rewrite \eqref{f2-final} in a different form. Observe that $\Psi(\eta)\cdot\eta=0$. Hence the numerator in \eqref{f2-final} can be written as follows:
\begin{equation}
\begin{split}
\Psi(\eta)+(\Psi(\eta)\cdot e_1)(u-e_1)
&=\Psi(\eta)+\frac{\Psi(\eta)\cdot x}{|\eta \times x|}\left(\sin\mu\eta-(1-\cos\mu) \frac{x-(\eta\cdot x)\eta}{|\eta \times x|}\right)\\
&=\Psi(\eta)+\frac{\Psi(\eta)\cdot x}{L}\left(\eta-\frac{1-\cos\mu}{\sin^2\mu} \frac{x-(\eta\cdot x)\eta}{L}\right)\\
&=\Psi(\eta)+\frac{\Psi(\eta)\cdot x}{L^2(1+\cos\mu)}((1+L)\eta-x)\\
&=\Psi(\eta)+\frac{\Psi(\eta)\cdot x}{L(1+\cos\mu)}\left(\eta+\frac{\eta-x}{|\eta-x|}\right).
\end{split}
\label{f2-int1}
\end{equation} 
Here we have used again that $\Psi(\eta)\cdot\eta=0$. Now the smooth dependence on $x$ in a neigborhood of the origin is obvious.


\section{General domains. Outline of argument.}\label{general_dom}

Let $D$ denote the convex domain where $\fbf$ is supported. The domain is supposed to be known.

It is easy to see that
\begin{equation}\label{rot-f1}  \nabla\times \fbf=\nabla\times \fbf^{(1)}.   \end{equation}
Indeed, by construction
\begin{equation}
\fbf(x)-\fbf^{(1)}(x)=\int_{S^2} \pa_p^2 \Rnor{\bf f}(p,\eta)|_{p=\eta\cdot x} d\eta=\int_{S^2} \eta \psi(\eta\cdot x,\eta)d\eta
\label{diff}
\end{equation} 
for some scalar function $\psi$.  Direct calculation shows that calculating the curl of the integral in \eqref{diff} is zero.

Observe that $\fbf^{(1)}(x)$ can be represented in the form
\begin{equation}
\fbf^{(1)}(x)=\int_{S^2} [\Psi(\eta\cdot x,\eta)-\eta(\eta\cdot \Psi(\eta\cdot x,\eta))]d\eta
\label{f1}
\end{equation} 
for some vector function $\Psi$. Direct calculation shows that $\nabla\cdot$ of the integral in \eqref{f1} is zero, i.e. $\nabla\cdot \fbf^{(1)}=0$.
 
Consequently, $\fbf^s(x)-\fbf^{(1)}(x)$ is a harmonic vector field, i.e. $\fbf^s(x)-\fbf^{(1)}(x)=\nabla h(x)$ for some $h$ such that $\Delta h\equiv 0$.  Here $\fbf^s$ is the solenoidal part of $\fbf$, and we used \eqref{diff} and that $\nabla\times \fbf=\nabla\times \fbf^s$.
 
Once $\fbf^{(1)}(x)$, $x\in D$, is computed, we can compute its cone beam transform and {\bf subtract} from the data. Since potential vector fields are in the kernel of the cone beam transform, we can think that the measured data is the cone beam transform of $\fbf^s$, not of $\fbf$. Thus the subtraction gives the cone beam transform of the harmonic vector field $\fbf^s(x)-\fbf^{(1)}(x)=\nabla h$, which we denote $D_h(y(s),\Theta)$. Here $y(s)$ and $\Theta$ are the position of the source and the direction of the ray, respectively. Let $x_{in}(y(s),\Theta)$ and $x_{out}(y(s),\Theta)$ be the points where the ray determined by $y(s)$ and $\Theta$ enters the domain $D$ and exists the domain $D$, respectively. Obviously, $h(x_{out}(y(s),\Theta))-h(x_{in}(y(s),\Theta))=D_h(y(s),\Theta)$. Thus, we know the differences between the values of $h$ for many pairs of points $(x_{in},x_{out})$ on the boundary. If the collection of lines corresponding to our data $D_h(y(s),\Theta)$ is sufficiently rich (which is the case, for example, when the trajectory consists of two orthogonal circles), then we can find $h$ on all the boundary of $D$ up to a constant. Hence we can solve the following boundary value problem: $\Delta h=0$ in $D$, $h_{\partial D}=\text{known}$, and then set $\fbf^{(2)}=\nabla h$. Since we compute the gradient, the fact that boundary values of $h$ are known only up to a constant does not affect the computation of $\fbf^{(2)}$.


\section{Numerical experiments}\label{numerics}

\noindent We present some implementations of the improved inversion formula (\ref{recon-3}), (\ref{f2-final}). We confine ourselves to smooth vector fields $\fbf$ supported in the closed unit ball $\overline{B^3}$. If a solenoidal vector field $\fbf$  vanishes at the boundary in the sense of  $\fbf(\eta)\cdot \eta=0$ for all $\eta \in \mathbb{S}^2$, then if follows that $\fbf^{(2)}= 0$, see also (\ref{rot-f1}). For this reason, we assume that the second part of the inversion formula $\fbfzwei$ mainly contains information about boundary values. The numerical tests should emphasize this phenomenon as well as the exactness of the formula. We mention that in contrast to $\fbfeins$ we have to evaluate $\Psi (\eta)$ which is the most elaborate part of the inversion formula.

\begin{figure}[H]
\centering
\begin{minipage}[b]{0.4\textwidth}
\includegraphics[width=\textwidth]{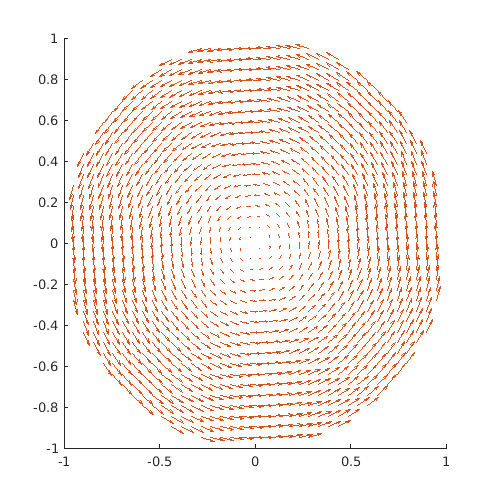}
\end{minipage}
\begin{minipage}[b]{0.4\textwidth}
\includegraphics[width=\textwidth]{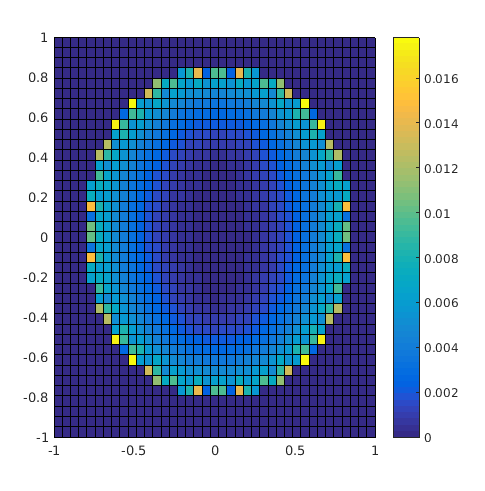}
\end{minipage}
\caption{The exact field $\fbf_a$ plotted in the plane $x_3= -0.5$ (left picture) and plot of $\|\fbfzwei\|$ for $x_3=-0.5$ (right picture).}
\label{fig-fa}
\end{figure}

The first vector field we reconstruct is given as
\[  \fbf_a (x) = \nabla \times \exp (-|x|^2/2) \begin{pmatrix} x_1\\x_2\\-x_3 \end{pmatrix}  \]
which is solenoidal and satisfies $\eta \cdot \fbf_a(\eta)=0$ on $S^2$. Hence we expect $\fbf_2$ to be zero. A plot of $\fbf_a$ for $x_3=-0.5$ can be seen in figure \ref{fig-fa} (left picture). The right picture in figure \ref{fig-fa} shows $\fbfzwei$ for this field and in fact demonstrates that this part of the inversion formula vanishes for $\fbf_a$ up to discretization errors.

%
%

\begin{figure}[H]
\centering
\includegraphics[scale=0.55]{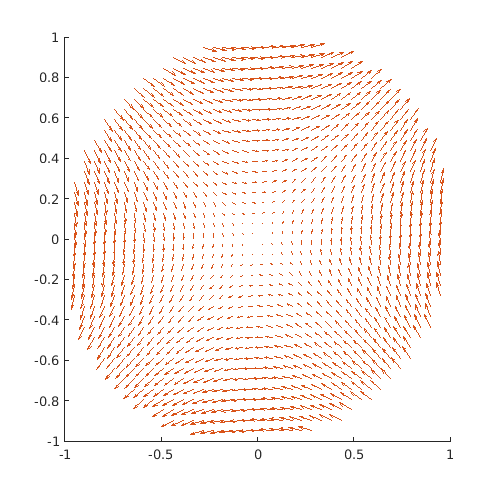}
\caption{The exact field $\fbf_b$ plotted in the plane $x_3=0$}
\label{fig-fb}
\end{figure}

The second vector field is given as
\[  \fbf_b (x) = \begin{pmatrix} x_2\\x_1\\0 \end{pmatrix}   \]
which is divergence free, either. A plot of $\fbf_b$ for $x_3=0$ is illustrated in figure \ref{fig-fb}.

\begin{figure}[H]
\centering
\begin{minipage}[b]{0.4\textwidth}
\includegraphics[width=\textwidth]{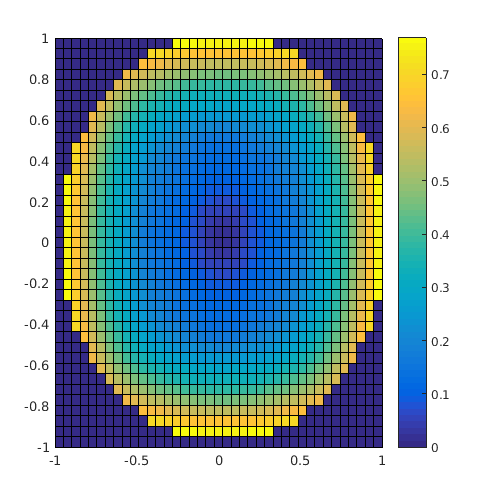}
\end{minipage}
\begin{minipage}[b]{0.4\textwidth}
\includegraphics[width=\textwidth]{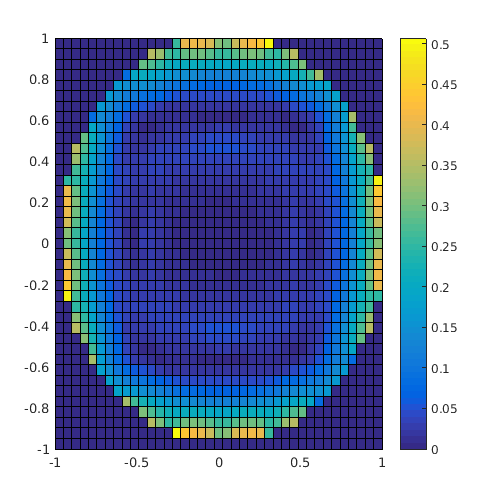}
\end{minipage}
\caption{Plot of $\|\fbfzwei\|$ in the plane $x_3=0$ (left picture) and absolute error $\|\fbf_b-(\fbfeins+\fbfzwei)\|$ in the plane $x_3=0$ (right picture).}
\label{fig-fb-f2}
\end{figure}

The visualization of $\fbfzwei$ in figure \ref{fig-fb-f2} for $\fbf_b$ in fact demonstrates that its main part is located close to the boundary of $B^3$.
the right picture in figure \ref{fig-fb-f2} shows the absolute error $\|\fbf_b-(\fbfeins+\fbfzwei)\|$ for $x_3=0$. Again we see that the biggest part of the error occurs at the boundary,
a phenomenon which was observed in other measure geometries, too, see e.g. \cite{pot_part}. The reasons for this are not entirely clarified. Besides discretization errors we think
that numerical instabilities occur in the integral of (\ref{f2-final}) for $x$ being close to the boundary $S^2=\partial B^3$. Numerical tests showed that the articfacts close to the boundary
appear also when computing $\fbf_1$.

\begin{figure}[H]
\centering
\includegraphics[scale=0.6]{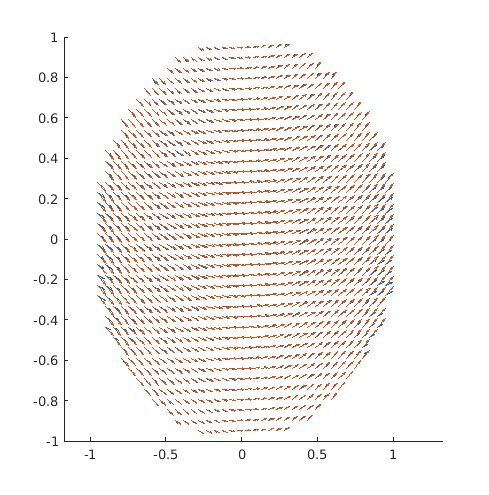}
\caption{Plot of $\fbf_c$ (red) versus $\fbfeins+\fbfzwei$ (blue) at $x_3=0$}
\label{fig-fc-inv}
\end{figure}

The third vector field is given by
\[ \fbf_c (x) = \begin{pmatrix} \cos(x_2)\\\sin(x_1)\\0 \end{pmatrix}\,. \]
Figure \ref{fig-fc-inv} shows $\fbf_c$ plotted versus $\fbfeins+\fbfzwei$ and demonstrates a high concurrence.

\begin{figure}[H]
\centering
\begin{minipage}[b]{0.4\textwidth}
\includegraphics[width=\textwidth]{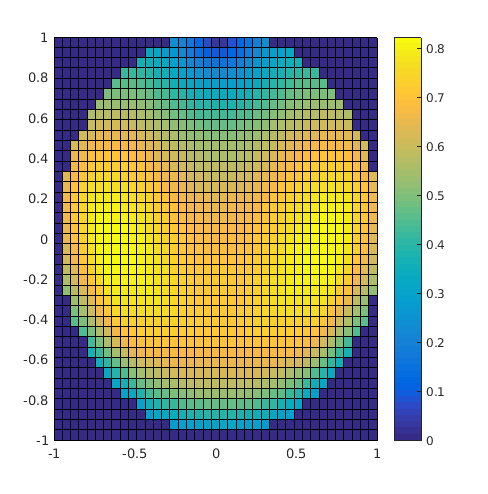}
\end{minipage}
\begin{minipage}[b]{0.4\textwidth}
\includegraphics[width=\textwidth]{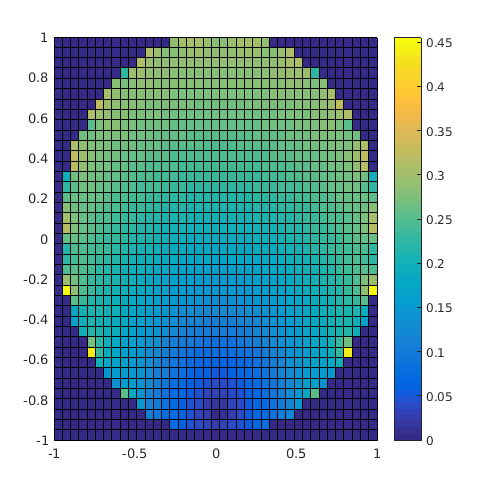}
\end{minipage}
\caption{Plot of $\|\fbfeins\|$ (left picture) and of $\|\fbfzwei\|$ (right picture) in the plane $x_3=0$}
\label{fig-fc-f1}
\end{figure}

\begin{figure}[H]
\centering
\includegraphics[scale=0.55]{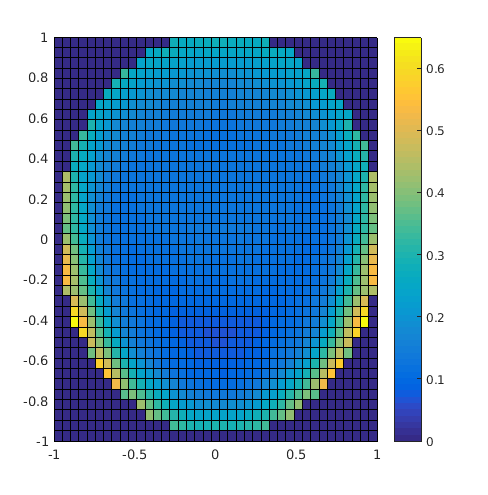}
\caption{Absolute error $\|\fbf_c-(\fbfeins+\fbfzwei)\|$ in the plane with $x_3=0$}
\label{fig-fc-error}
\end{figure}

Figure \ref{fig-fc-f1} presents plots of $\|\fbfeins\|$ and $\|\fbfzwei\|$, respectively, in the plane $x_3=0$. A look at these plots clearly demonstrates that $\fbfeins$ contains most information of the interior values of $\fbf_c$, whereas $\fbfzwei$ has its largest values close to the boundary. The absolute error $\|\fbf_c-(\fbfeins+\fbfzwei)\|$ (Figure \ref{fig-fc-error}) shows again high accuracy with discretization errors near $\partial B^3$.

The last vector field is given by
\[ \fbf_d (x) = \begin{pmatrix} x_2^2-x_3^2\\x_1^2-x_3^2\\0 \end{pmatrix}  \]

\begin{figure}[H]
\centering
\includegraphics[scale=0.6]{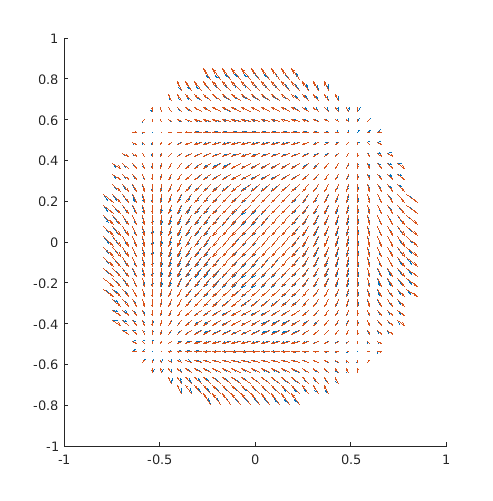}
\caption{Plot of $\fbf_d$ (red) versus $\fbfeins+\fbfzwei$ (blue) in the plane $x_3 = 0.5$}
\label{fig-fd}
\end{figure}

A plot of $\fbf_d$ versus the reconstruction $\fbfeins+\fbfzwei$ is shown in figure \ref{fig-fd}.

\begin{figure}[H]
\begin{minipage}[b]{0.4\textwidth}
\includegraphics[width=\textwidth]{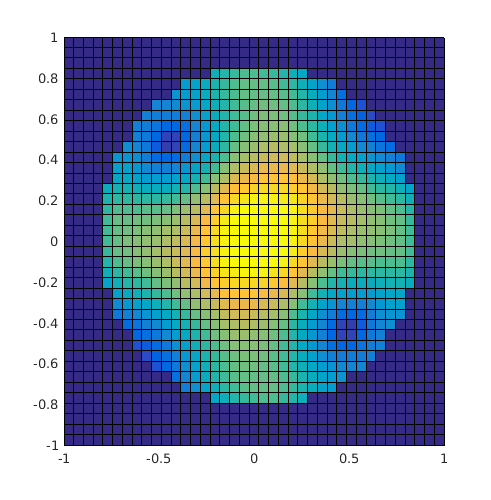}
\end{minipage}
\begin{minipage}[b]{0.4\textwidth}
\includegraphics[width=\textwidth]{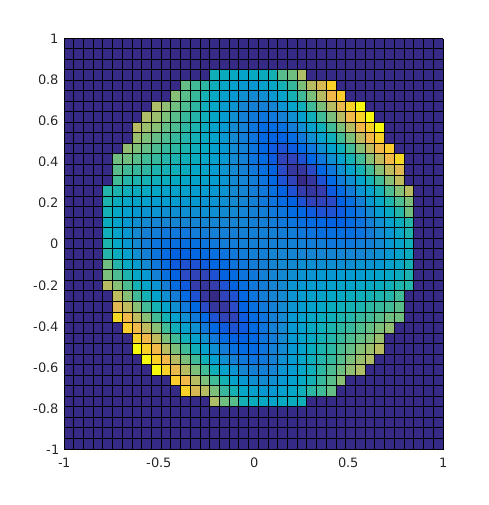}
\end{minipage}
\caption{Plot of $\|\fbfeins\|$ (left picture) and $\|\fbfzwei\|$ (right picture) in the plane $x_3=0.5$.}
\label{fig-fd-f1-f2}
\end{figure}

\begin{figure}[H]
\centering
\includegraphics[scale=0.55]{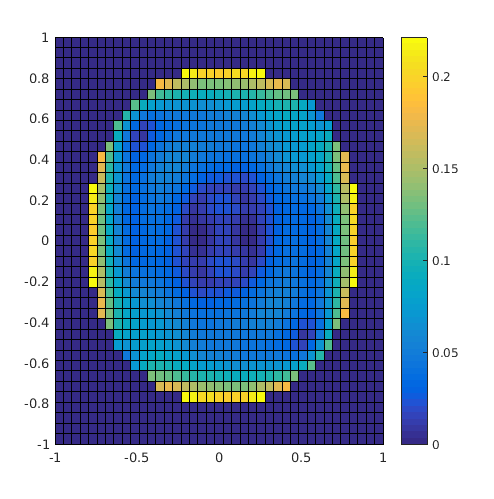}
\caption{Absolute error $\|\fbf_d-(\fbfeins+\fbfzwei)\|$ in the plane $x_3 = 0.5$}
\label{fd-error}
\end{figure}

Figure \ref{fig-fd-f1-f2} shows plots of $\|\fbfeins\|$ and $\|\fbfzwei\|$. These pictures again emphasize that $\fbfzwei$ mainly contributes to the boundary
values just as suggested by our theoretical investigations. Figure \ref{fd-error} finally illustrates the accuracy of our inversion formula up to
discretization errors. Again the values close to the boundary are most sensible with respect to errors. 

\section{Conclusions}

We improved the exact inversion formula for cone beam vector tomography achieved in \cite{katsch13} by saving one integration order
in the second part $\fbfzwei$ of this formula. Theoretical considerations suggest that the first part of the formula, $\fbfeins$,
which is of classical, filtered backprojection type, contains information abut the curl of $\fbf$, whereas the second part
$\fbf$ mainly contributes to the boundary values. Consequently $\fbfzwei=0$ for a solenoidal vector field with vanishing
boundary values $\fbf (\eta)\cdot\eta = 0$ on $S^2$. These theoretical investigations as well as a good performance of the
inversion formula were supported by numerical experiments
for different divergence free vector fields.
Future work will address general convex domains and the extension of the inversion formula to distributions.


\section*{Acknowledgements}

Alexander Katsevich and Thomas Schuster have been supported by German Science Foundation (Deutsche Forschungsgemeinschaft, DFG) under grant Schu 1978/12-1.

\bibliographystyle{amsalpha}
\bibliography{bibliogr,references-thomas}

\end{document}